\documentclass[12pt]{amsart}
\usepackage{amssymb}
\makeatletter
\@namedef{subjclassname@2020}{%
  \textup{2020} Mathematics Subject Classification}
\makeatother
\newtheorem{theorem}{Theorem}[section]
\newtheorem{lemma}[theorem]{Lemma}
\newtheorem{proposition}[theorem]{Proposition}
\theoremstyle{definition}
\newtheorem{definition}[theorem]{Definition}
\newtheorem{remark}[theorem]{Remark}
\newcommand{\CC}{\mathbb C}
\newcommand{\QQ}{\mathbb Q}
\newcommand{\RR}{\mathbb R}
\newcommand{\ZZ}{\mathbb Z}
\newcommand*{\set}[1]{\{#1\}}
\newcommand{\OO}{\mathbb O}
\usepackage{xcolor} 
\usepackage{url}
\usepackage[colorlinks=true,linkcolor=blue,citecolor=purple]{hyperref}
\begin{document}
\baselineskip=17pt
\title[Euler's magic matrices]{On Euler's magic matrices of sizes $3$ and $8$}
\author[P. Müller]{Peter M\"uller}
\address{Institute of Mathematics\\ University of W\"urzburg\\ Germany}
\email{peter.mueller@uni-wuerzburg.de}
\date{}
\begin{abstract}
  A proper Euler's magic matrix is an integer $n\times n$ matrix
  $M\in\ZZ^{n\times n}$ such that $M\cdot M^t=\gamma\cdot I$ for some
  nonzero constant $\gamma$, the sum of the squares of the entries
  along each of the two main diagonals equals $\gamma$, and the
  squares of all entries in $M$ are pairwise distinct. Euler
  constructed such matrices for $n=4$. In this work, we use
  multiplication matrices of the octonions to construct
  examples for $n=8$, and prove that no such matrix exists for $n=3$.
\end{abstract}
\maketitle
\section{Introduction} A classical magic square is an $n\times n$
matrix $A$ with distinct nonnegative integer entries such that the
sums of the entries in each row, each column, and both main diagonals
are the same.

If one requires in addition that the entries of $A$ are square
numbers, then $A$ is called a \emph{magic square of squares}, see
\cite{boyer}, \cite{robertson}. It is an open problem whether a magic
square of squares of size $3$ exists, despite considerable effort on
this question in \cite{bremner_1} and \cite{bremner_2}.

Leonhard Euler looked at the question for $n=4$. He noticed that
orthogonal matrices, or slightly more general matrices $M$ such that
$M\cdot M^t=\gamma\cdot I$, could make the problem easier. For if $A$
is the matrix whose entries are the squares of those of $M$, then the
conditions on the row and column sums for $A$ are automatically
fulfilled. So one is faced with only two polynomial conditions for the
two main diagonals, and the added requirement that the entries of $A$
are pairwise distinct.
\begin{definition}\label{def:emm}
  Let $R$ be a commutative ring and $n$ be a positive integer. A
  matrix $M=(m_{i,j})_{1\le i, j\le n}\in R^{n\times n}$ is called an
  \emph{Euler's magic matrix} over $R$, if the following holds for
  some $\gamma\in R\setminus\set{0}$, where $I_n$ denotes the
  $n\times n$ identity matrix:
  \begin{align}
    \label{eq:1}
    M\cdot M^t &= \gamma\cdot I_n,\\
    \sum_{i=1}^nm_{i,i}^2 &= \gamma,\label{eq:2}\\
    \sum_{i=1}^nm_{i,n+1-i}^2 &= \gamma.\label{eq:3}
  \end{align}
  If in addition the squares of the entries in $M$ are pairwise
  distinct, then we call $M$ a \emph{proper Euler's magic matrix}.
\end{definition}
Note that the sum in \eqref{eq:3} is the sum of the squares of the
elements on the anti-diagonal. As remarked above, if $M$ is a proper
Euler's magic matrix, then $A=(a_{i,j})_{1\le i,j\le n}$ with
$a_{i,j}=m_{i,j}^2$ is a magic square of squares.

The case $n=2$ is easy to work out, the only Euler's magic matrices
have the form $\begin{pmatrix}a&a\\a&-a\end{pmatrix}$,
$\begin{pmatrix}a&a\\-a&a\end{pmatrix}$,
$\begin{pmatrix}a&-a\\a&a\end{pmatrix}$, and
$\begin{pmatrix}-a&a\\a&a\end{pmatrix}$. In particular, there are no
proper such cases.

Euler studied the intriguing problem of whether a proper Euler's magic
matrix exists for $n=4$, referring to it as a ``problema
curiosum''. In fact he managed to produce a multi-parametric family of
such matrices. A particular case, due to Euler too, is\[
  M = \begin{pmatrix}68& -29& 41& -37\\
    -17& 31& 79& 32\\
    59& 28& -23& 61\\
    -11& -77& 8& 49
  \end{pmatrix}.
\]
A more systematic study for the case $n=4$, based on the algebra of
quaternions, was given by Hurwitz in \cite[Vorlesung
12]{hurwitz:quaternionen}. See \cite{boyer}, \cite{pirsic:msb} and in
particular \cite{oswald_steuding} and \cite[Lecture
12]{oswald_steuding:book} about the history of this problem and its
connection to the Hamilton quaternions.

In Section 2 we show this negative result concerning $3\times3$ matrices:
\begin{theorem}\label{t:3x3}
  There is no Euler's magic matrix in $\QQ^{3\times 3}$.
\end{theorem}
It is perhaps surprising that $n=3$ is the only positive integer for
which there is no Euler's magic matrix in $\QQ^{n\times n}$. In fact,
in Section \ref{n=5} we provide a simple construction of Euler's magic
matrices in $\ZZ^{n\times n}$ for each $n\ne3$. However, these
matrices are far from proper. In particular, we do not know whether
a proper Euler's magic matrix exists for $n=5$. Section \ref{n=5}
contains some examples which come close.

In Section 3 we examine the case $n=8$. In \cite{pirsic:arxiv} and
\cite{pirsic:msb}, \'Isabel Pirsic suggested to use certain matrices
coming from left and right multiplication of the octonions to find a
construction analogous to Euler's for $n=4$. However, a massive search
by her did not find an example. In this note, we demonstrate that a
careful analysis of the polynomial system yields solutions without
much searching. So Pirsic's suggestion to use octonion multiplication
matrices proved successful.

In order to describe our results (and later the methods), we introduce
the matrices
\[L(a,b,\dots,h)=
\left(\begin{array}{rrrrrrrr}
a & -b & -c & -d & -e & -f & -g & -h \\
b & a & -d & c & -f & e & h & -g \\
c & d & a & -b & -g & -h & e & f \\
d & -c & b & a & -h & g & -f & e \\
e & f & g & h & a & -b & -c & -d \\
f & -e & h & -g & b & a & d & -c \\
g & -h & -e & f & c & -d & a & b \\
h & g & -f & -e & d & c & -b & a
\end{array}\right)
\]
and
\[R(p,q,\dots,w)=
\left(\begin{array}{rrrrrrrr}
p & -q & -r & -s & -t & -u & -v & -w \\
q & p & s & -r & u & -t & -w & v \\
r & -s & p & q & v & w & -t & -u \\
s & r & -q & p & w & -v & u & -t \\
t & -u & -v & -w & p & q & r & s \\
u & t & -w & v & -q & p & -s & r \\
v & w & t & -u & -r & s & p & -q \\
w & -v & u & t & -s & -r & q & p
\end{array}\right).
\]
With respect to a suitable basis of the octonions $\OO$ over the real
numbers, $L(a,b,\dots,h)$ describes the left multiplication
$\OO\to\OO, z\mapsto xz$, where $x$ has the coefficients
$a,b,\dots,h$. Likewise, $R(p,q,\dots,w)$ describes the right
multiplication $\OO\to\OO, z\mapsto zx$, where $x$ has coefficients
$p,q,\dots,w$.

The reason that we use different letters for the entries of $L(\cdot)$
and $R(\cdot)$ is that the matrices we will study have the form
$M=L\cdot R=L(a,b,\dots,h)\cdot R(p,q,\dots,w)$. The point is that for
$L$ and $R$ (with arbitrary $a,b,\dots,h,p,q,\dots,w$) we have
$L\cdot L^t=(a^2+b^2+\dots+h^2)\cdot I_8$ and
$R\cdot R^t=(p^2+q^2+\dots+w^2)\cdot I_8$, and therefore
$M\cdot M^t=\gamma\cdot I_8$ with
$\gamma=(a^2+b^2+\dots+h^2)(p^2+q^2+\dots+w^2)$.

Thus condition \eqref{eq:1} in Definition \ref{def:emm} is
automatically satisfied, and one has ``only'' to discuss the two
polynomial conditions \eqref{eq:2} and \eqref{eq:3} in the $16$
unknowns $a,b,\dots,h,p,q,\dots,w$ and the properness. A typical
example that we obtain is
\begin{theorem}
  Set 
  \begin{align*}
    L &= L(2,\; 1,\; 1,\; 4,\; 2,\; 1,\; 1,\; -2),\\
    R &= R(-7,\; -55,\; -11,\; 1,\; -27,\; -13,\; -19,\; 4).
  \end{align*}
  Then
  \[
    M = L\cdot R=%
\left(\begin{array}{rrrrrrrr}
142 & 197 & -225 & 30 & 16 & 57 & -13 & -170 \\
-37 & -60 & 136 & 201 & 177 & 98 & -32 & -193 \\
-283 & -4 & -148 & -95 & 71 & 164 & 10 & -1 \\
-22 & 237 & 181 & -178 & 138 & -29 & -45 & -8 \\
-120 & 97 & 27 & 74 & -62 & -129 & 293 & -82 \\
-9 & 38 & -116 & 131 & 235 & -144 & 0 & 187 \\
-103 & 180 & 50 & 195 & -163 & 64 & -132 & 107 \\
-126 & -35 & -51 & -20 & -42 & -247 & -191 & -154
\end{array}\right)
\]
is a proper Euler's magic matrix.
\end{theorem}
This single example arises from specializing the parameters of a
$4$-parametric family of Euler's magic matrices to values which
preserve properness. In this case, the example was obtained from
setting $(q, r, t, u) = (-55, -11, -27, -148)$ and rescaling in the
following theorem.
\begin{theorem}\label{t:8x8}
  For variables $q, r, t, u$ over $\QQ$ set
\[
  X = 7 q^{2} + 7 r^{2} + 21 q t - 7 r t + 34 t^{2} - 7 q u - 21 t%
  u + 4 u^{2} + 7 q + 21 r - 7 u + 34
\]
  and
  \begin{align*}
    L &= L(2,\; 1,\; 1,\; 4,\; 2,\; 1,\; 1,\; -2),\\
    R &= R\left(\frac{3(t^2 - 1)u}{2X},\; q,\; r,\; 1,\; t,\;%
        u-q-3t-1,\; t - r - 3,\; \frac{u^2-X}{2u}\right).
  \end{align*}
  Then $M=L\cdot R$ is a proper Euler's magic matrix over
  $\QQ(q,r,t,u)$.
\end{theorem}
\begin{remark}
  The reader who wishes to verify the examples need not type these
  matrices. A proof is provided in the ancillary SageMath script
  \texttt{euler\_verify.sage} at \cite{euler_verify}. This script can
  be run at the SageMathCell at \url{https://sagecell.sagemath.org/}.
\end{remark}
\section{There are no Euler's magic \texorpdfstring{$3\times3$}{3}
  matrices}
If we look for Euler's magic $n\times n$ matrices for odd $n$ over a
field, then we may assume that these matrices are orthogonal:
\begin{lemma}\label{lem:ortho}
  Let $K$ be a field, $n$ be odd, and $M\in K^{n\times n}$ with
  $M\cdot M^t=\gamma\cdot I_n$ for $0\ne\gamma\in K$. Then
  $\gamma=\lambda^2$ for $\lambda\in K$, therefore
  $\frac{1}{\lambda}\cdot M$ is orthogonal.
\end{lemma}
\begin{proof}
  Write $n=2k+1$. Taking the determinant in
  $M\cdot M^t=\gamma\cdot I_n$ yields
  $(\det M)^2=\gamma^n=\gamma^{2k+1}$, so $\gamma=\lambda^2$ for
  $\lambda=\det M/\gamma^k$.
\end{proof}
We use the following refinement of the Cayley transform which
parametrizes orthogonal real matrices:
\begin{proposition}\label{prop:cayley}
  Let $M\in\RR^{n\times n}$ be a an orthogonal matrix. Then one can write
  \[
    DM=(I_n-S)(I_n+S)^{-1},
  \]
  where $S\in\RR^{n\times n}$ is skew-symmetric and
  $D\in\RR^{n\times n}$ is a diagonal matrix with diagonal entries in
  $\set{-1,1}$.
\end{proposition}
Note that $0$ is the only possible real eigenvalue of a real
skew-symmetric matrix, hence $I_n+S$ is invertible for every real
skew-symmetric matrix $S$.
\begin{remark}
  This proposition was stated and proved (in a slightly different
  form) in 1991 by Liebeck and Osborne
  \cite{liebeck_osborbe:amm}. However, it had already appeared 30
  years earlier in \cite[Chapter 6, Section 4, Exercises 7 --
  11]{bellman:book}. We briefly sketch the argument: The Cayley
  transform $S\mapsto M=(I_n-S)(I_n+S)^{-1}$ is a bijection from the
  set of skew-symmetric matrices $S\in\RR^{n\times n}$ to the set of
  orthogonal matrices $M\in\RR^{n\times n}$ for which $-1$ is not an
  eigenvalue. This map is involutive in the sense that if
  $M=(I_n-S)(I_n+S)^{-1}$, then $S=(I_n-M)(I_n+M)^{-1}$. All of this
  has been well known since Cayley's time (and is easy to verify). To
  prove the proposition, one needs to find $D$ such that $-1$ is not
  an eigenvalue of $DM$. Since $DM+I_n=D\cdot(M+D)$, this if
  equivalent to $M+D$ being invertible. But the existence of $D$
  follows from an easy induction on $n$ (for arbitrary
  $M\in\RR^{n\times n}$), see, e.g., \cite[Lemma 1]{hsu:matrices}.
\end{remark}

We now prove Theorem \ref{t:3x3}. Let $M\in\QQ^{3\times3}$ be an
Euler's magic matrix. Multiplying $M$ by a nonzero rational preserves
this property, so in view of Lemma \ref{lem:ortho} we may assume that
$M$ is orthogonal. Moreover, the property of $M$ being an Euler's
magic matrix is preserved upon replacing rows by their
negatives. Thus, by Proposition \ref{prop:cayley}, we may assume that
\[
  M = (I_3-S)(I_3+S)^{-1},\text{ where }S=\begin{pmatrix}0&a&b\\
  -a&0&c\\
    -b&-c&0\end{pmatrix}\in\QQ^{3\times3}.
\]
With $\Delta=\det(I_3+S)=a^2+b^2+c^2+1$ we compute
\[
  M=\frac{1}{\Delta}\left(\begin{array}{rrr}
    -a^{2} - b^{2} + c^{2} + 1 & -2 b c - 2 a & 2 a c - 2 b \\
    -2 b c + 2 a & -a^{2} + b^{2} - c^{2} + 1 & -2 a b - 2 c \\
    2 a c + 2 b & -2 a b + 2 c & a^{2} - b^{2} - c^{2} + 1
  \end{array}\right).
\]
The conditions \eqref{eq:2} and \eqref{eq:3} about the diagonal and
the anti-diagonal are $D=E=0$ with
\begin{align*}
  D %
  &= (-a^2 - b^2 + c^2 + 1)^2%
    + (-a^2 + b^2 - c^2 + 1)^2%
    + (a^2 - b^2 - c^2 + 1)^2 - \Delta^2\\
  &= 2(a^4 - 2a^2b^2 + b^4 - 2a^2c^2 - 2b^2c^2 + c^4 - 2a^2 - 2b^2 - 2c^2 + 1)
\intertext{and}
  E %
  &= (2ac - 2b)^2%
    + (-a^2 + b^2 - c^2 + 1)^2%
    + (2ac + 2b)^2-\Delta^2\\
  &= 4(-a^2b^2 + 2a^2c^2 - b^2c^2 - a^2 + 2b^2 - c^2).
\end{align*}
Now, while $D$ and $E$ look somewhat complicated, we get by some
magical calculation
\[
  \frac{D+E}{2} = (a^2 - 2b^2 + c^2 - 2)^2-3(b^2 + 1)^2.
\]
From $D=E=0$ and the fact that $3$ is not a square in $\QQ$ we get
$b^2+1=0$, a contradiction.
\begin{remark}
  Reducing to orthogonal matrices and using the Cayley transform is
  quite natural and straightforward. But how could one have guessed
  that $D=E=0$ with $D$ and $E$ as above has no rational solution? We
  sketch another, less ``magical'', proof.

  The first thing one might observe is that $D$ and $E$ are polynomials
  in $a^2$, $b^2$, and $c^2$, and that $D$ and $E$ are symmetric in
  $a^2$ and $c^2$. Thus we can express $D$ and $E$ in terms of
  $\beta=b^2$, $s=a^2+c^2$, and $p=a^2c^2$. We obtain
  \begin{align*}
    D/2 &= \beta^2 -2(1+s)\beta + (1-s)^2 - 4p\\
    E/4 &= (2-s)\beta - s + 2p.
  \end{align*}
Eliminating $\beta$ from $D=E=0$ yields
\begin{align*}
  0 &= 4 p^{2} + \left(-8 s^{2} + 16 s - 8\right) p + s^{4} - 4 s^{3}
      + 12 s^{2} - 16 s + 4\\
  &= 4(p-(s-1)^2)^2 -3(s - 2)^2s^2.
\end{align*}
As $3$ is not a square in $\QQ$, we get $s=0$ or $s=2$, and
then $p=(s-1)^2=1$. The case $s=0$ yields $a^2+c^2=0$. Then $a=c=0$,
and therefore $p=a^2c^2=0$, a contradiction.

In the other case we have $a^2+c^2=s=2$ and $a^2c^2=p=1$,
hence
\[
  (a^2-1)^2+(c^2-1)^2=(a^2+c^2)^2%
  -2(a^2+c^2)+2-2a^2c^2=0.
\]
We get $a^2=c^2=1$. This yields $D/4=b^4 - 6b^2 - 3=(b^2-3)^2-12$, but
there is no $b\in\QQ$ with $D=0$.
\end{remark}
\section{ Euler's magic \texorpdfstring{$8\times8$}{8} matrices}
We look for $16$ integers $a,b,\dots,h,p,q,\dots,w$ such that
$M=L\cdot R$ with $L=L(a,b,\dots,h)$ and $R=R(p,q,\dots,w)$ is a
proper Euler's magic matrix.

As remarked in the introduction, we have
\[
  M\cdot M^t=(a^2+b^2+\dots+h^2)(p^2+q^2+\dots+w^2)I_8,
\]
so \eqref{eq:1} in Definition \ref{def:emm} holds with
$\gamma=(a^2+b^2+\dots+h^2)(p^2+q^2+\dots+w^2)$. Thus the conditions
\eqref{eq:2} and \eqref{eq:3} are to be studied. In the $3\times 3$
case from the previous section, it turned out that actually the sum of
the equations \eqref{eq:2} and \eqref{eq:3} had a useful property. In
the situation here it again appears that the sum and the difference of
\eqref{eq:2} and \eqref{eq:3} have somewhat better properties.

Accordingly, with $M=(m_{i,j})_{1\le i,j\le 8}$, set
\begin{align*}
  A(a,b,\dots,w) &= \sum_{i=1}^8m_{i,i}^2-\sum_{i=1}^8m_{i,9-i}^2,\\
  B(a,b,\dots,w) &= \sum_{i=1}^8m_{i,i}^2+\sum_{i=1}^8m_{i,9-i}^2-2\gamma.
\end{align*}
Thus $M$ is an Euler's magic matrix if and only if $A(a,b,\dots,w)=0$
and $B(a,b,\dots,w)=0$.
\subsection{\texorpdfstring{Some properties of $A(a,b,\dots,w)$ and
    $B(a,b,\dots,w)$}{A and B}}

The strategy is the following. We fix integers $a,b,\dots,h$ and
consider $p,q,\dots,w$ as variables. Then each entry of $M=(m_{i,j})$
is a linear form in $p,q,\dots,w$, and therefore $A$ and $B$ are
homogeneous quadratic forms.

To ease the language, we call an arbitrary matrix \emph{proper} if and
only if the squares of its entries are pairwise distinct.

Also, it is obvious that if a matrix which depends on parameters is
not proper, then this is even more true if we specialize parameters.

So in order to start with, we need to choose $a,b,\dots,h\in\ZZ$ such
that $M\in\ZZ[p,q,\dots,w]^{8\times 8}$ is proper.

This for instance requires that $a\ne0$ or $h\ne0$, as
$m_{1,8}-m_{8,1}=2(aw-hp)$. In fact if two of the numbers
$a,b,\dots,h$ vanish, then $M$ is not proper, as one can easily check
with a simple program.

But certain other choices of $a,b,\dots,h$ are seen to be impossible
even if $M$ is proper. For instance $M$ is proper for
$a=b=c=d=e=f=g=h=1$. But in this case we get
\[
  A(p,q,\dots,w) = (p + q + t + u)(r + s + v + w),
\]
which forces $p + q + t + u=0$ or $r + s + v + w$. However
\begin{align*}
  m_{3,3}-m_{3,6} &= 2(p + q + t + u),\\
  m_{2,2}-m_{2,7} &= 2(r + s + v + w),
\end{align*}
so no matter which factor vanishes, we see that the condition
$A(p,q,\dots,w)=0$ forces $M$ to be improper. The same happens (for
other index pairs of $M$) whenever all $a,b,\dots,h$ are in
$\set{-1,1}$.

However, if at least one of the integers $a,b,\dots,h$ is different
from $\pm1$ (and assuming without loss of generality that
$\gcd(a,b,\dots,h)=1$), it rarely happens that the quadratic form $A$
(or $B$) is reducible.

A necessary condition for $A=B=0$ of course is that the quadratic
forms $A$ and $B$ are isotropic, and in fact that every linear
combination $\lambda A+\mu B$ is isotropic. However, that is a very
weak condition, because we have more than $4$ variables, so the only
condition is that $\lambda A+\mu B$ is isotropic over $\RR$. But for
random choices of $a,b,\dots,h$ this is usually the case.
\subsection{A naive search}
If we eliminate one of the variables, $w$ say, from $A=B=0$, then in
general we will obtain a quartic form in $p,q,\dots,v$. Next, if we
specialize all but two of these variables to integers, then usually
the resulting curve in the remaining variables will be quartic. Many
experiments have shown that the rational points of this quartic are
hard to analyze. Despite having degree $4$, its genus (if it is
absolutely irreducible) will be at most $1$. (It is a known fact from
algebraic geometry that if the intersection of two quadratic surfaces
in $\CC^3$ is an irreducible curve, then its genus is at most
$1$. See, e.g., \cite[Lecture 22, Pencils of
Quadrics]{harris:ag_first_course}.) However, in general the curve will
have no rational singularities which could help to transform it to a
cubic. But even in cases when there were rational singularities, and
furthermore the transformed cubic could be transformed to Weierstrass
normal form, present state software like SageMath, Pari, or Magma was
not able to compute the Mordell-Weil rank of these curves in all cases
we tried. The easily computed torsion points led to no examples. Note
that even if we have a rational point on the curve, then $w$ usually
has degree $2$ over $\QQ$, because in general $A$ has $w$-degree
$2$. And in most cases where we found rational solutions of $A=B=0$,
the resulting matrix $M$ was not proper.

In rare cases however we found valid solutions by this approach. For
instance for
\begin{align*}
  (a,b,\dots,h) &= (0, 1, 1, 1, 1, 1, -1, 5)\\
  (p,q,r,s,t) &= (3, -2, -4, 5, 6)
\end{align*}
the resulting system $A(u,v,w)=B(u,v,w)=0$ has the rational solution
$(u,v,w)=(13/15, -14/15, -23/5)$, which indeed gives, after rescaling,
the proper Euler's magic square
\[
  L(0, 1, 1, 1, 1, 1, -1, 5)\cdot R(45, -30, -60, 75, 90, 13, -14,
  -69).
\]
Finding a few examples like this required checking thousands of
potential integer tuples $(a,b,\dots,h,p,q,r,s,t)$ of length $13$.
Here we used a combination of a backtracking and a greedy algorithm to
find tuples $(a,b,\dots,h,p,q,r,s,t)$ of small integers such that the
specialized matrix $M\in\QQ[u,v,w]^{8\times 8}$ is still proper. For
if the integers are large, it is less likely that the system $A=B=0$
has rational solutions $(u, v, w)$.
\subsection{Too strong restrictions}
As there are so many more variables than equations, one might consider
to impose strong restrictions to make the polynomial system more
manageable. For instance, one might pick $3$ distinct variables
$X,Y,Z\in\set{p,q,\dots,w}$, and hope to specialize the remaining
$13$ variables to integers such that $A(X, Y, Z)$, as a polynomial in
$X, Y, Z$, is $0$.

For instance, the coefficients in $A$ of $w^2$ and $wp$ are
$8(h-a)(h+a)$ and $16ah$, respectively. So if want them both to
vanish, then $a=h=0$, in which case $M$ is not proper anymore.

So we cannot have $\set{p, w}\subset\set{X,Y,Z}$. Many other combinations
fail for the analogous reason. But some cases require a finer
analysis.
\subsection{A working restriction}
The following compromise proved to be fruitful. We look for conditions
on $a,b,\dots,h\in\ZZ$ such that $A$ and $B$ both have degree $1$ in
$w$. Assume this for a moment. Let $x$ and $y$ be the coefficients of
$w$ in $A$ and $B$, respectively. Then $F=yA-xB$ is a cubic form
$F\in\ZZ[p,q,\dots,v]$, where we have eliminated $w$.

Besides having lower degree, the added advantage is that a solution of
$F=0$ with $p,q,\dots,v\in\QQ$ extends to a solution of $A=B=0$,
provided that $x,y\ne0$.

Fortunately, the condition for $A$ and $B$ to have $w$-degree $\le1$ is
rather easy and not very restrictive: As noted above, the coefficient
of $w^2$ in $A$ is $8(h-a)(h+a)$. So we need to pick $h=\pm
a$. Furthermore, the coefficient of $w^2$ in $B$ is, up to the factor
$-2$, equal to
$b^2 + c^2 + d^2 + e^2 + f^2 + g^2-3(a^2+h^2)=b^2 + c^2 + d^2 + e^2 +
f^2 + g^2-6a^2$. Thus $A$ and $B$ have $w$-degree $\le1$ if and only
if $h=\pm a$ and $b^2 + c^2 + d^2 + e^2 + f^2 + g^2 = 6a^2$.

As remarked previously, if $a=h=0$, then $M$ is not proper. Thus we
assume that
\begin{equation}\label{eq:w1}
  h = \pm a\ne0,\;\;\; b^2 + c^2 + d^2 + e^2 + f^2 + g^2 = 6a^2.
\end{equation}
Next let $F=yA-xB$ be the cubic form in $\ZZ[p,q,\dots,v]$. The idea
is to find a rational specialization of some of the variables
$p,q,\dots,v$ such that $F$ will have degree $1$ with respect to one
of these variables. Because then we get an immediate rational
parametrization of $F=0$ and hence of the solutions of $A=B=0$.

It happens that \eqref{eq:w1} already implies that $F$ has degree at
most $2$ in $p$, because the coefficient of $p^3$ in $F$ is
$32h\cdot(b^2 + c^2 + d^2 + e^2 + f^2 + g^2 - 6a^2)=0$.

Also, the coefficient of $p^2$ in $F$ has, up to the nonzero factor
$-128h^2$, the useful form
\begin{align*}
  \left(a g + b h\right) q + \left(-a f + c h\right) r + \left(-a e +
  d h\right) s\\
  + \left(a d + e h\right) t + \left(a c + f h\right) u + \left(-a b + g h\right) v.
\end{align*}
As $h=\pm a\ne0$ and not both $b$ and $g$ vanish, we get that
$ag+bh\ne0$ or $-ab+gh\ne0$. So we can solve for $q$ or $v$ in terms
of the remaining variables to achieve that $F$ has degree
at most $1$ in $p$. Usually this degree equals $1$, so we can solve
for $p$ to make $F$ vanish. Finally, provided that $A$ and $B$ after these
substitutions still have degree $1$ in $w$, solving for $w$ finally
yields $A=B=0$, where now the entries of the matrix $M$ are rational
functions over $\QQ$ in the variables $r, s, t, u$ and $q$ or $v$.

This is essentially how we got the parametrization in Theorem
\ref{t:8x8}, starting with $(a,b,\dots,h)=(2, 1, 1, 4, 2, 1, 1, -2)$.
\begin{remark}
  The integer tuples $(a, b,\dots, h)$ satisfying \eqref{eq:w1} can
  easily be parametrized. So if we use such a parametrization instead
  of the specific tuple like $(2, 1, 1, 4, 2, 1, 1, -2)$ which led to
  Theorem \ref{t:8x8}, we can still carry out the procedure. If we
  consider matrices over $\QQ$ as equivalent if they differ by a
  nonzero scalar factor, then we obtain an $8$-parametric family of
  proper Euler's magic matrices for which Theorem \ref{t:8x8} is just
  a subcase.
\end{remark}
\section{\texorpdfstring{Improper Euler's $n\times n$ magic matrices for $n\ge4$}{improper}}\label{n=5}
Let $\sigma$ be an element of the symmetric group on
$\set{1,2,\dots,n}$. We define the matrix
$M(\sigma)=(m_{i,j})\in\ZZ^{n\times n}$ by
\[
  m_{i,j}=
  \begin{cases}
    1,& \text{if } j=\sigma(i)\\
    0,& \text{otherwise.}
  \end{cases}
\]
Note that $M(\sigma)$ contains exactly one $1$ in each row and column,
hence $M(\sigma)\cdot M(\sigma)^t=I_n$.
\begin{theorem}
  Set \[ \sigma=
    \begin{cases}
      (1\;2\;\ldots\;n-1)(n),& \text{if }n\text{ is }even\\
      (1\;2\;\ldots\;k-1)(k)(k+1\;k+2\;\ldots n),& \text{if }n=2k-1%
                                                \text{ is odd.}\\
    \end{cases}
  \]
Then $M(\sigma)$ is an Euler's magic matrix.
\end{theorem}
\begin{proof}
  We need to verify the conditions \eqref{eq:1}, \eqref{eq:2}, and
  \eqref{eq:3}. As $M(\sigma)\cdot M(\sigma)^t=I_n$, condition
  \eqref{eq:1} holds with $\gamma=1$.
  
  First suppose that $n=2k\ge4$ is even. Then $m_{i,i}=0$ for
  $1\le i<n$ and $m_{n,n}=1$. Furthermore, $m_{i,n+1-i}=1$ if and only
  if $i=k$. We see that $M(\sigma)$ fulfills conditions \eqref{eq:2}
  and \eqref{eq:3}.

  Similarly, if $n=2k-1$ is odd, then $m_{k,k}=1$ and $m_{i,i}=0$ if
  $i\ne k$, and $m_{i, n+1-i}=1$ if and only if $i=k$. Again,
  $M(\sigma)$ fulfills conditions \eqref{eq:2} and \eqref{eq:3}.
\end{proof}
Of course, these matrices are far from being proper. We have tried to
find proper Euler's magic matrices for $n=5$. The approach was to use
the Cayley transform as in the case $n=3$. This results in two
polynomial conditions in $10$ variables over $\QQ$. Or, after
homogenization, we need to solve two polynomials over the integers in
$11$ unknowns. As all variables appear in a high degree, we basically
tried a random search. Surprisingly, this way we obtained quite a few
Euler's magic matrices, and some of them came close to properness. In
fact, the following are Euler's magic matrices, where the squares of
its entries give $24$ distinct elements. In each case we highlight the
pair of entries which violates properness:
\[
  \left(\begin{array}{rrrrr}
    -106 & -32 & -8 & -75 & -50 \\
    -4 & -38 & -120 & 58 & -35 \\
    24 & \mathbf{20} & -73 & -88 & 80 \\
    61 & 66 & -16 & -46 & -100 \\
    70 & -115 & \mathbf{20} & -40 & -18
  \end{array}\right)
\]
\[
  \left(\begin{array}{rrrrr}
4 & 3 & 40 & -94 & -142 \\
-29 & -128 & -90 & 44 & -58 \\
154 & 28 & -35 & 56 & -42 \\
74 & \mathbf{-82} & -10 & -114 & 73 \\
-24 & \mathbf{82} & -140 & -61 & 2
  \end{array}\right)
\]
\[
  \left(\begin{array}{rrrrr}
    -204 & -38 & 10 & -11 & -312 \\
    54 & -262 & -260 & 36 & -13 \\
    -84 & \mathbf{102} & -165 & -306 & 48 \\
    291 & \mathbf{102} & -40 & -56 & -202 \\
    66 & -223 & 210 & -206 & -2
  \end{array}\right)
\]
\[
  \left(\begin{array}{rrrrr}
    29 & -218 & -370 & \mathbf{-188} & 180 \\
    88 & -384 & 22 & 158 & -269 \\
    -160 & 58 & -40 & -333 & -334 \\
    210 & -139 & 304 & -286 & 124 \\
    418 & \mathbf{188} & -147 & -4 & -146
  \end{array}\right)
\]
\[
  \left(\begin{array}{rrrrr}
    \mathbf{-392} & -336 & -21 & 282 & 210 \\
    177 & -384 & -24 & \mathbf{-392} & 240 \\
    408 & -186 & -246 & 357 & -40 \\
    -48 & -309 & 176 & -42 & -510 \\
    -192 & 14 & -546 & -168 & -165
  \end{array}\right)
\]
A verification of these examples is again provided in \cite{euler_verify}.
\subsection*{Acknowledgement}I thank the anonymous referee for helpful
comments, including pointing out errors and suggestions that
improved the manuscript.

\end{document}